\documentclass[a4paper,12pt]{amsart}
\usepackage{amsfonts}

\usepackage{amssymb}
\usepackage{ifthen}
\usepackage{graphicx}
\usepackage{tablists}
\nonstopmode \numberwithin{equation}{section}
\usepackage[usenames]{color}
\usepackage{enumerate}
\newtheorem{thm}[equation]{Theorem}
\newtheorem{lem}[equation]{Lemma}
\newtheorem{cor}[equation]{Corollary}

\newtheorem{cl}{Claim}[section]
\newtheorem{ca}{Case}[section]
\newtheorem{sca}{Subcase}[section]
\newtheorem{scl}[section]{Subclaim}
\newtheorem{conj}[equation]{Conjecture}

\theoremstyle{definition}
\newtheorem{defn}[equation]{Definition}
\newtheorem{prob}[equation]{Problem}
\newtheorem{op}[equation]{Open Problem}
\newtheorem{ques}[equation]{Question}
\newtheorem{rem}[equation]{Remark}
\newtheorem{exam}[equation]{Example}

\newcounter {own}
\def\theown {\thesection       .\arabic{own}}

\makeatletter
\@namedef{subjclassname@2020}{\textup{2020} Mathematics Subject Classification}
\makeatother

\newenvironment{pf}[1][]{%
 \vskip 3mm
 \noindent
 \ifthenelse{\equal{#1}{}}%
  {{\slshape Proof. }}%
  {{\slshape #1.} }%
 }%
{\qed\bigskip}

\newcommand{\id}{{\operatorname{id}}}

\def\be{\begin{equation}}
\def\ee{\end{equation}}

\newcommand{\ben}{\begin{enumerate}}
\newcommand{\een}{\end{enumerate}}

\newcommand{\blem}{\begin{lem}}
\newcommand{\elem}{\end{lem}}
\newcommand{\bthm}{\begin{thm}}
\newcommand{\ethm}{\end{thm}}
\newcommand{\bcor}{\begin{cor}}
\newcommand{\ecor}{\end{cor}}
\newcommand{\beg}{\begin{exam}}
\newcommand{\eeg}{\end{exam}}
\newcommand{\begs}{\begin{examples}}
\newcommand{\eegs}{\end{examples}}
\newcommand{\bdefe}{\begin{defn}}
\newcommand{\edefe}{\end{defn}}
\newcommand{\bprob}{\begin{prob}}
\newcommand{\eprob}{\end{prob}}
\newcommand{\bques}{\begin{ques}}
\newcommand{\eques}{\end{ques}}
\newcommand{\bei}{\begin{itemize}}
\newcommand{\eei}{\end{itemize}}
\newcommand{\bcon}{\begin{conj}}
\newcommand{\econ}{\end{conj}}
\newcommand{\bop}{\begin{op}}
\newcommand{\eop}{\end{op}}

\newcommand{\bca}{\begin{ca}}
\newcommand{\eca}{\end{ca}}
\newcommand{\bsca}{\begin{sca}}
\newcommand{\esca}{\end{sca}}

\newcommand{\bcl}{\begin{cl}}
\newcommand{\ecl}{\end{cl}}

\newcommand{\bscl}{\begin{scl}}
\newcommand{\escl}{\end{scl}}

\newcommand{\bcons}{\begin{conjs}}
\newcommand{\econs}{\end{conjs}}
\newcommand{\bprop}{\begin{propo}}
\newcommand{\eprop}{\end{propo}}
\newcommand{\br}{\begin{rem}}
\newcommand{\er}{\end{rem}}
\newcommand{\brs}{\begin{rems}}
\newcommand{\ers}{\end{rems}}
\newcommand{\bo}{\begin{obser}}
\newcommand{\eo}{\end{obser}}
\newcommand{\bos}{\begin{obsers}}
\newcommand{\eos}{\end{obsers}}
\newcommand{\bpf}{\begin{pf}}
\newcommand{\epf}{\end{pf}}
\newcommand{\ba}{\begin{array}}
\newcommand{\ea}{\end{array}}
\newcommand{\beq}{\begin{eqnarray}}
\newcommand{\beqq}{\begin{eqnarray*}}
\newcommand{\eeq}{\end{eqnarray}}
\newcommand{\eeqq}{\end{eqnarray*}}

\newcommand{\ra}{\rightarrow}

\newcounter{minutes}
\divide\time by 60
\newcounter{hours}
\multiply\time by 60 \addtocounter{minutes}{-\time}

\begin{document}

\bibliographystyle{amsplain}
\title{Gromov hyperbolization of unbounded spaces and Hamenst\"adt metric}

\author{Qingshan Zhou}
\address{Qingshan Zhou, School of Mathematics and Big Data, Foshan University,  Foshan, Guangdong 528000, People's Republic
of China} \email{qszhou1989@163.com; q476308142@qq.com}

\author{Saminathan Ponnusamy}
\address{Saminathan Ponnusamy, Department of Mathematics, Indian Institute of Technology Madras, Chennai 600036,
India}
\email{samy@iitm.ac.in}

\author{Qianghua Luo}
\address{Qianghua Luo, School of Mathematics and Big Data, Foshan University,  Foshan, Guangdong 528000, People's Republic
of China} \email{15616203413@163.com}

\def\thefootnote{}
\footnotetext{ \texttt{\tiny File:~\jobname .tex,
          printed: \number\year-\number\month-\number\day,
          \thehours.\ifnum\theminutes<10{0}\fi\theminutes}
} \makeatletter\def\thefootnote{\@arabic\c@footnote}\makeatother

\date{}

\begin{abstract} In this paper, we investigate Gromov hyperbolizations of unbounded locally complete and  incomplete metric spaces associated with three hyperbolic type metrics: the hyperbolization metric introduced by Ibragimov, the distance ratio metric, and the quasihyperbolic metric. As an application, we obtain a Gromov hyperbolic characterization of unbounded uniform domains in Banach spaces.
\end{abstract}

\thanks{Qingshan Zhou was supported by by Guangdong Basic and Applied Basic Research Foundation (No. 2021A1515012289)}

\maketitle{} \pagestyle{myheadings} \markboth{Zhou et al.}{Gromov hyperbolization of unbounded spaces and Hamenst\"adt metric}

\section{Introduction and main results}\label{sec-1}
The theory of $\delta$-hyperbolic spaces, described in \cite{Gr87}, is based on an important observation that the asymptotic properties of classical hyperbolic space $\mathbb{H}^n$ can be ensured by using a simple condition for quadruples of points. Ignoring the local structure, Gromov introduced an inequality between arbitrary four points in a metric space and built up a theory of general negative curvature type spaces, now known as Gromov hyperbolic spaces. Many classical results concerning the large scale geometry of simply connected complete manifolds of negative curvature (such as $\mathbb{H}^n$) can be generalized to this class of hyperbolic spaces. See \cite{BS,BrHa,BuSc,D10,FS11,Gr87,Z1,SX,X14,Vai} for more backgrounds and motivations on this area.

It is known that every simply connected planar proper sub-domain admits a Gromov hyperbolic metric by pulling back the classical hyperbolic metric of the unit disk via the Riemann mapping. Recently, Ibragimov introduced a new metric that hyperbolizes (in the sense of Gromov) all locally compact incomplete metric spaces in \cite{Z1}. Subsequently, it is shown that the Gromov boundary $\partial_h X$ of $(X,h)$ can be identified with the metric boundary $\partial X$ of a bounded locally compact incomplete metric space $(X,d)$ via a quasisymmetric map; for the unbounded case, the extended boundary $\partial X\cup\{\infty\}$ is equipped with a chordal metric, see \cite[Theorem 3.1]{Z2}.

As in \cite{Z2}, by a Gromov hyperbolization of a locally complete and incomplete metric space $(X,d)$ we mean the following: Equipping $X$ with a metric $m$ such that the space $(X,m)$ is Gromov hyperbolic and the Gromov boundary $\partial_m X$ of $(X,m)$ can be identified with the metric boundary $\partial X$ of $(X,d)$ via a quasisymmetric map.

To complement \cite[Theorem 3.1]{Z2}, we focus our attention on the problem how to hyperbolize unbounded incomplete metric spaces.  Note that the local compactness property of the space is \textbf{not} assumed throughout this paper.
It is known that quasisymmetric mappings map unbounded sets to unbounded sets. However, the Gromov boundary of a Gromov hyperbolic space equipped with a common visual metric is bounded, see \cite[p.\,435]{BrHa} or \cite[Section 6]{BS}.

To overcome this limitation, we introduce the class of Hamenst\"adt metrics on the punctured Gromov boundary of a hyperbolic space which was considered by Hamenst\"adt \cite{Ha} in studying negatively pinched Hadamard manifolds and defined via horospherical distances. Contrary to the case of visual metrics, the punctured Gromov boundary of a hyperbolic space $X$ equipped with a Hamenst\"adt metric is unbounded. It is very similar to the Poincar\'{e} half--space model of $\mathbb{H}^n$, see \cite{BuSc}.

The class of Hamenst\"adt metrics are defined by using Busemann functions, for related definitions and properties see \cite[Section 3.3]{BuSc}. Roughly speaking, a Busemann function on a Gromov hyperbolic space is defined to be the distance function from a point on
the Gromov boundary. This notion is very useful in many areas. For instance, employing Busemann functions on CAT$(-1)$ spaces, Foertsch and Radke \cite{FR11} characterized complete CAT($\kappa$) spaces with $\kappa<0$, with geodesic Hamenst\"adt boundary up to isometry. Moreover, Foertsch and  Schroeder \cite{FS11} investigated the relationship between Gromov hyperbolic spaces, CAT(-1) spaces and the Ptolemy inequality on their Gromov boundaries.

One of them is called a parabolic visual metric based on the vertical geodesic in some negatively curved solvable Lie groups in \cite{SX}. In \cite{SX, X14}, Shanmugalingam and Xie proved that all self quasi-isometries of these groups are almost isometries. It should be noted that this parabolic visual metric was formerly named Euclid--Cygan metric by Hersonsky and Paulin \cite{HP} in the study of the rigidity of discrete isometry groups of negatively curved spaces. With the aid of this notion, Dymarz \cite{D10,D14} recently studied the quasi-isometric rigidity of mixed type locally compact amenable hyperbolic groups.

As our first main result, we show that the hyperbolization metric $h$ (see (\ref{h})), introduced by Ibragimov in \cite{Z1}, hyperbolizes all unbounded locally complete and incomplete metric spaces as follows.

\begin{thm}\label{thm-1}
Let $(X,d)$ be an unbounded locally complete and incomplete metric space. Then there is a natural $\eta$-quasisymmetric map
$$\varphi\colon (\partial_h X\setminus\{\xi\},h_{b,\varepsilon})\to (\partial X,d),$$
where $\partial_h X\setminus\{\xi\}$ is the punctured Gromov boundary of hyperbolic space $(X,h)$ equipped with a Hamenst\"adt metric $h_{b,\varepsilon}$ and $b\in \mathcal{B}_h(\xi)$ is a Busemann function with  $\varphi(\xi)=\infty$. The function $\eta$ depends only on $\varepsilon$.
\end{thm}

\br
We recall the definition of natural maps (cf. \cite{BHK,HSX,Vai}). Let $(X,d)$ be a locally complete and incomplete metric space and let $m$ be a metric on $X$ such that $(X,m)$ is Gromov hyperbolic, and let $\partial_m X$ be the Gromov boundary of $(X,m)$. Let $\hat{X}$ be the one point extension $X\cup\{\infty\}$ of $X$ if $X$ is unbounded, and $ \hat{X}=X$ if $X$ is bounded. Then $\partial  \hat{X}=\partial X$ if $X$ is bounded, and $\partial  \hat{X}=\partial X\cup\{\infty\}$ if $X$ is unbounded.  Following \cite{Vai}, the topology of $\hat{X}$ consists of all open sets in $X$ and of all sets $U$ containing $\infty$ such that $\hat{X}\backslash U=X\backslash U$ is closed and bounded in $X$.

If the identity map $\id\colon (X,m)\to (X,d)$ has a continuous extension $\psi\colon \partial_m X\to \partial  \hat{X}$, then $\psi$ is called a {\it natural} map. Moreover, if such a natural map $\psi$ exists,  then we say that $\partial_m X$ is naturally equivalent to  $\partial  \hat{X}$.
\er

In addition to the classical hyperbolic metric, there are many other Gromov hyperbolic type metrics defined and extensively used in geometric function theory; for instance, the quasihyperbolic metric $k$ (see (\ref{k})) and the distance ratio metric $\widetilde{j}$ which were introduced by Gehring and Osgood in \cite{GO}. In \cite{Ha03}, H\"ast\"o demonstrated that any proper domain in $\mathbb{R}^n$ endowed with the $\widetilde{j}$-metric (see (\ref{j})) is $\delta$-hyperbolic. As  an application of Theorem \ref{thm-1}, we  prove that $\widetilde{j}$-metric also hyperbolizes all incomplete unbounded metric spaces as follows.

\begin{thm}\label{c-1} Let $(X,d)$ be an unbounded locally complete and incomplete metric space. Then  $(X, \widetilde{j})$ is Gromov hyperbolic and there is a natural $\widetilde{\eta}$-quasisymmetric map
$$\widetilde{\varphi}\colon (\partial_{\widetilde{j}} X\setminus\{\xi\},\widetilde{j}_{b,\varepsilon})\to (\partial X,d),$$
where $\partial_{\widetilde{j}} X\setminus\{\xi\}$ is the punctured Gromov boundary of hyperbolic space $(X,{\widetilde{j}})$ equipped with a Hamenst\"adt metric ${\widetilde{j}}_{b,\varepsilon}$ and $b\in \mathcal{B}_{\widetilde{j}}(\xi)$ is a Busemann function with $\widetilde{\varphi}(\xi)=\infty$. The function $\widetilde{\eta}$ depends only on $\varepsilon$.
\end{thm}

\br
Note that our argument in establishing the Gromov hyperbolicity of $(X,{\widetilde{j}})$ is different with \cite{Ha03}. Recently, the authors in \cite{ZhPG22} proved a result similar to Theorem \ref{c-1} in Euclidean spaces setting. Our approach is also different with that of \cite{ZhPG22}. According to \cite[Theorem 3.1]{Z1}, the identity map $\id\colon  (X,{\widetilde{j}})\to (X,h)$ is roughly similar. We demonstrate that roughly similar maps preserve the Gromov hyperbolicity of metric spaces (not necessarily geodesic), which induces a quasisymmetric correspondence  between the punctured Gromov boundaries equipped with Hamenst\"adt metrics, see Lemma \ref{zz-01}. Therefore, Theorem \ref{c-1} follows immediately from Theorem \ref{thm-1} and the above facts.
\er

It was shown in \cite[Theorem 3.6]{BHK} by Bonk et al. that a locally compact bounded uniform metric space is Gromov hyperbolic in the quasihyperbolic metric and the metric boundary is naturally quasisymmetrically equivalent to the Gromov boundary.
The second goal of this paper is to obtain an analogue of \cite[Theorem 3.6]{BHK} and to hyperbolize unbounded uniform metric spaces associated  with the quasihyperbolic metric. Our result is the following:

\begin{thm}\label{thm-2}
Let $A\geq 1$ and let $(X,d)$ be an unbounded $A$-uniform metric space. Then there is a natural $\eta$-quasisymmetric map
$$\phi\colon (\partial_k X\setminus\{\xi\},k_{b,\varepsilon})\to (\partial X,d),$$
where $\partial_k X\setminus\{\xi\}$ is the punctured Gromov boundary of hyperbolic space $(X,k)$ equipped with a Hamenst\"adt metric $k_{b,\varepsilon}$ and $b\in \mathcal{B}_k(\xi)$ is a Busemann function with  $\phi(\xi)=\infty$. The function $\eta$ depends only on $A$ and $\varepsilon$.
\end{thm}

In \cite{BHK}, Bonk et al. further studied the relationship between Gromov hyperbolic domains and uniform domains of $\mathbb{R}^n$.  Recall that a proper domain in $\mathbb{R}^n$ is called Gromov hyperbolic if it is $\delta$-hyperbolic with respect to its quasihyperbolic metric for some $\delta\geq 0$.  They showed that a bounded domain in $\mathbb{R}^n$ is uniform if and only if it is both Gromov hyperbolic and its Euclidean boundary is naturally quasisymmetrically equivalent to the Gromov boundary, see \cite[Theorem 1.11]{BHK}. Subsequently, V\"ais\"al\"a \cite{Vai} generalized this result to Banach spaces and obtained a dimensional free result. Instead of bounded domains he considered also unbounded domains, where the  quasisymmetric equivalence is replaced by quasim\"obius  equivalence.

In this paper, with the aid of Theorem \ref{thm-2}, we establish an unbounded analogue of \cite[Theorem 1.11]{BHK} and prove that an unbounded domain in a Banach space is uniform if and only if it is both Gromov hyperbolic and its norm boundary is naturally quasisymmetrically equivalent to the punctured Gromov boundary equipped with a Hamenst\"adt metric.

\begin{thm}\label{c-2} Let $D$ be an unbounded proper domain $($an open and connected set$)$ in a Banach space $(E,|\cdot|)$ with dimension $\dim E\geq 2$. Then $D$ is $A$-uniform if and only if $(D,k)$ is $\delta$-hyperbolic and there is a natural $\eta$-quasisymmetric map
$$\phi\colon (\partial_k D\setminus\{\xi\},k_{b,\varepsilon})\to (\partial D,|\cdot|),$$
where $\partial_k D\setminus\{\xi\}$ is the punctured Gromov boundary of hyperbolic space $(D,k)$ equipped with a Hamenst\"adt metric $k_{b,\varepsilon}$ and $b\in \mathcal{B}_k(\xi)$ is a Busemann function with $\phi(\xi)=\infty$. The parameters $A$ and $\delta$, $\eta$, $\varepsilon$ depend only on each other.
\end{thm}

The rest of this paper is organized as follows.  In Section \ref{sec-2}, we recall some definitions and preliminary results. The proofs of Theorems \ref{thm-1} and \ref{c-1} are given in Section \ref{sec-3}. Finally, our goal in Section \ref{sec-4} is to show Theorems \ref{thm-2} and  \ref{c-2}.

\section{Preliminaries and auxiliary results}\label{sec-2}
\subsection{Notions and notations}
Following \cite{BuSc}, for $t_1, t_2, c\in \mathbb{R}$ with $c\geq 0$, it is convenient to write $t_1\doteq t_2$ up to an error $\leq c$ or $t_1\doteq_c t_2$ instead of $|t_1-t_2|\leq c$.

In what follows, $(X,d)$ denotes a metric space with the metric $d$. We often write the distance between $x$ and $y$ as $d(x,y)$, and $d(x,A)$ to denote the distance from a point $x$ to a set $A$. Also, we denote the open metric ball with center $x\in X$ and radius $r>0$ by $B(x,r)=\{z\in X:\; d(z,x)<r\}$.

Let $\overline{X}$ be the metric completion of a metric space $X$ and let $\partial X= \overline{X}\setminus X$ denote the metric boundary of $X$. If $\partial X\neq \emptyset$, then  $X$ is {\it incomplete} and we denote $d(x)=d(x,\partial X)$ for all $x\in X$. Note that $d(x)$ may vanish such as  the case of $X=\mathbb{Q}\subset \mathbb{R}$, where $\mathbb{Q}$ denotes the set of all rational numbers in $\mathbb{R}$. To avoid this from happening, we assume that the space $(X,d)$ is {\it locally complete}; that is, every point in $X$ has a complete neighborhood. Thus for a locally complete and incomplete metric space $(X,d)$, $d(x)>0$ for all $x\in X$.  Obviously, locally compact spaces are locally complete.

By a curve, we mean a continuous function $\gamma\colon $ $[s,t]\to X$. For a curve $\alpha$, we let $\alpha[u,v]$ be  a fixed but arbitrary subcurve of $\alpha$ between two points $u,$ $v\in \alpha$. The length of $\gamma$ is denoted by
$$\ell(\gamma)=\sup\left\{\sum_{i=1}^{n}d(\gamma(t_i),\gamma(t_{i-1}))\right\},$$
where the supremum is taken over all finite partitions $s=t_0<t_1<t_2<\cdots<t_n=t$. The curve $\gamma$ is called {\it rectifiable} if $\ell(\gamma)<\infty$. Also, $X$ is called {\it rectifiably connected} if each pair of points in the space can be joined by a rectifiable curve.

The length function $s_{\gamma}$, associated with a rectifiable curve $\gamma\colon [t_1,t_2]\to X$, is  defined by
$s_{\gamma}(t)=\ell(\gamma|_{[t_1,t]})$ for $t\in [t_1,t_2]$. For any rectifiable curve $\gamma\colon [t_1,t_2]\to X$, there is a unique parametrization $\gamma_s\colon $ $[0, \ell(\gamma)]\to X$ such that $\gamma=\gamma_s\circ s_{\gamma}$. Obviously, $\ell(\gamma_s|_{[0,t]})=t$ for $t\in [0, \ell(\gamma)]$. The parametrization $\gamma_s$ is called the {\it arclength parametrization} of $\gamma$.

The {\it cross ratio} of  a quadruple $Q = (x,y,z,w)$ of four distinct points  in $X$ is defined to be the number
$$\tau(Q) = |x,y,z,w| = \frac{d(x,z)d(y,w)}{d(x,y)d(z,w)}.$$

Observe that the definition is extended in the well known manner to the case
where one of the points is $\infty$. For example,
$$|x,y,z, \infty| = \frac{d(x,z)}{d(x,y)}.$$

\bdefe A homeomorphism $f\colon $ $(X,d)\to (X',d')$ between metric spaces is called
\begin{enumerate}
\item $\eta$-{\it quasisymmetric} if there is a homeomorphism $\eta \colon  [0,\infty) \to [0,\infty)$ such that
$$ d(x,y)\leq td(x,z)\;\; \mbox{implies}\;\;   d'(f(x),f(y)) \leq \eta(t)d'(f(x),f(z))$$
for each $t>0$ and for each triple $x,$ $y$, $z$ of points in $X$;

\item {\it $\theta$-quasim\"obius} if there is a homeomorphism $\theta \colon  [0,\infty) \to [0,\infty)$ such that
$$\tau(f(Q))\leq\theta(\tau(Q))$$ holds for each quadruple $Q\subset X$.
\end{enumerate}\edefe
Note that the properties of quasisymmetric and quasim\"obius maps can be found in \cite{TV, Vai-5, Vai-4}.

\subsection{Gromov hyperbolic geometry}
In this subsection, we give some basic information about Gromov hyperbolic spaces, see e.g. \cite{BrHa,BuSc,Gr87,Vai-0}. Let $(X,d)$ be a metric space. We say that $X$ is {\it Gromov hyperbolic} if there is a constant $\delta\geq 0$ such that
$$(x|y)_w\geq \min\{(x|z)_w,(z|y)_w\}-\delta$$
for all $x,y,z,w\in X$, where $(x|y)_w$ is the {\it Gromov product} with respect to $w$ defined by
$$(x|y)_w=\frac{1}{2}[d(x,w)+d(y,w)-d(x,y)].$$

\bdefe Suppose that $(X, d)$ is a Gromov hyperbolic space and $w\in X$.\ben
\item
A sequence $\{x_i\}$ in $X$ is called a {\it Gromov sequence} if $(x_i|x_j)_w\rightarrow \infty$ as $i,$ $j\rightarrow \infty$.
\item
Two such Gromov sequences $\{x_i\}$ and $\{y_i\}$ are said to be {\it equivalent} if $(x_i|y_i)_w\rightarrow \infty$  as $i\rightarrow \infty$.
\item
The {\it Gromov boundary} $\partial_\infty X$ of $X$ is defined to be the set of all equivalence classes of Gromov sequences.
\item
For $x\in X$ and $\xi\in \partial_\infty X$, the Gromov product $(x|\xi)_w$ of $x$ and $\xi$ is defined by
$$(x|\xi)_w= \inf \big\{ \liminf_{i\rightarrow \infty}(x|y_i)_w:\; \{y_i\}\in \xi\big\}.$$
\item
For $\xi,$ $\zeta\in \partial_\infty X$, the Gromov product $(\xi|\zeta)_w$ of $\xi$ and $\zeta$ is defined by
$$(\xi|\zeta)_w= \inf \Big\{ \liminf_{i\rightarrow \infty}(x_i|y_i)_w:\; \{x_i\}\in \xi\;\;{\rm and}\;\; \{y_i\}\in \zeta\Big\}.$$
\een
\edefe

\begin{lem}\label{z00}$($\cite[Lemma $5.11$]{Vai-0}$)$
Let $X$ be a $\delta$-hyperbolic space, let $o,z\in X$, and let $\xi,\xi'\in\partial_\infty X$. Then for any sequences $\{y_i\}\in \xi$, $\{y_i'\}\in \xi'$, we have
\begin{enumerate}
\item  $(z|\xi)_o\leq  \liminf_{i\rightarrow \infty} (z|y_i)_o \leq  \limsup_{i\rightarrow \infty} (z|y_i)_o\leq (z|\xi)_o+\delta;$
\item  $(\xi|\xi')_o\leq  \liminf_{i\rightarrow \infty} (y_i|y_i')_o \leq  \limsup_{i\rightarrow \infty} (y_i|y_i')_o\leq (\xi|\xi')_o+2\delta.$
\end{enumerate}
\end{lem}

For $0<\varepsilon<\min\{1,\frac{1}{5\delta}\}$, we define
$$\rho_{w,\varepsilon}(\xi,\zeta) =e^{-\varepsilon(\xi|\zeta)_w}$$
for all $\xi,\zeta$ in the Gromov boundary $\partial_\infty X$ of $X$ with convention $e^{-\infty}=0$. We now define
$$d_{w,\varepsilon}(\xi,\zeta):=\inf\bigg\{\sum_{i=1}^{n} \rho_{w,\varepsilon} (\xi_{i-1},\xi_i):n\geq 1,\xi=\xi_0,\xi_1,\ldots,\xi_n=\zeta\in \partial_\infty X\bigg\}.
$$
Then $(\partial_\infty X,d_{w,\varepsilon})$ is a metric space with
\be\label{t-3} \frac{1}{2}\rho_{w,\varepsilon}\leq d_{w,\varepsilon} \leq \rho_{w,\varepsilon},\ee
and we call $d_{w,\varepsilon}$ a {\it visual metric} of $\partial_\infty X$ with base point $w\in X$ and parameter $\varepsilon>0$.

Next, we introduce the class of Hamenst\"adt metrics based at a Busemann function or a point on the Gromov boundary of hyperbolic spaces. Following \cite{BuSc}, we say that $b\colon X\to \mathbb{R}$ is a {\it Busemann function} based at  $\xi\in \partial_{\infty}X$, denoted by $b=b_{\xi,w}\in \mathcal{B}(\xi)$ for some $w\in X$, if for all $x\in X$, we have
$$b(x)=b_{\xi,w}(x)=(\xi|w)_x-(\xi|x)_w.$$

Moreover, we define the Gromov product of $x,y\in X$ based at the Busemann function $b=b_{\xi,w}\in \mathcal{B}(\xi)$ by
$$(x|y)_b=\frac{1}{2}\big(b(x)+b(y)-d(x,y)\big).$$
Similarly, for $x\in X$ and $\zeta\in \partial_\infty X$, the Gromov product $(x|\zeta)_b$ of $x$ and $\zeta$ is defined by
$$(x|\zeta)_b= \inf \Big\{ \liminf_{i\rightarrow \infty}(x|z_i)_b:\; \{z_i\}\in \zeta\Big\}.
$$
For points $(\xi_1,\xi_2)\in (\partial_\infty X\times  \partial_\infty X)  \setminus (\xi,{\xi})$, we define their Gromov product based at $b$ via the formula
$$(\xi_1|\xi_2)_b=\inf\big\{\liminf_{i\to\infty} (x_i|y_i)_b: \{x_i\}\in\xi_1 , \{y_i\}\in\xi_2\}.$$

According to \cite[$(3.2)$ and Example $3.2.1$]{BuSc}, we see that
\be\label{z0} (x|y)_b\doteq_{10\delta} (x|y)_{\xi,w}=(x|y)_w-(\xi|x)_w-(\xi|y)_w.\ee
It follows from \cite[Proporsition $3.2.3$]{BuSc} that $(x|y)_b,(y|z)_b,(x|z)_b$ form a $22\delta$-triple for all $x,y,z\in X$.  Let $C\geq 0$. By a $C$-triple we mean that a triple of real numbers with the property that the two smallest of these numbers differ by at most $C$.

In view of these facts, we recall the definition of {\it Hamenst\"adt metrics} of  $\partial_\infty X\setminus\{\xi\}$ which is based at $\xi$ or a Busemann function $b=b_{\xi,w}\in \mathcal{B}(\xi)$. For $\varepsilon>0$ with $e^{22\varepsilon\delta}\leq 2$, we define
$$\rho_{b,\varepsilon}(\xi_1,\xi_2)= e^{-\varepsilon(\xi_1|\xi_2)_b}\;\;\;\;\;\;\;\mbox{for all}\;\xi_1,\xi_2\in \partial_\infty X\setminus\{\xi\}.$$
Then for $i=1,2,3$ with $\xi_i\in \partial_\infty X\setminus\{\xi\}$, we have
$$\rho_{b,\varepsilon}(\xi_1,\xi_2)\leq e^{22\varepsilon\delta} \max\{\rho_{b,\varepsilon}(\xi_1,\xi_3),\rho_{b,\varepsilon}(\xi_3,\xi_2)\}.$$
Let $\xi_1$, $\xi_2\in\partial_\infty X\setminus\{\xi\}$, we now define
$$d_{b,\varepsilon}(\xi_1,\xi_2):=\inf\bigg\{\sum_{i=1}^{n} \rho_{b,\varepsilon} (\zeta_{i-1},\zeta_i):n\geq 1,\xi_1=\zeta_0,\zeta_1,...,\zeta_n=\xi_2\in\partial_\infty X\setminus\{\xi\}\bigg\}.
$$
Again by \cite[Lemma $3.3.3$]{BuSc}, it follows that $(\partial_\infty X\setminus\{\xi\}, d_{b,\varepsilon})$ is a metric space with
\be\label{z-0.1} \frac{1}{2}\rho_{b,\varepsilon} \leq d_{b,\varepsilon}\leq \rho_{b,\varepsilon}.\ee
Then $d_{b,\varepsilon}$ is called a {\it Hamenst\"adt metric} on the punctured space $\partial_\infty X\setminus\{\xi\}$ based at $b$ with parameter $\varepsilon$.

\subsection{Quasihyperbolic metric and uniform spaces} In this part, we assume that $(X,d)$ is  an incomplete  and rectifiably connected metric space. The {\it quasihyperbolic length} of a rectifiable curve $\gamma$ in $X$ is given by
$$\ell_{k}(\gamma)=\int_{\gamma}\frac{ds}{d(z)}.
$$
For any $x,y$ in $X$,
the {\it quasihyperbolic distance}
$k(x,y)$  between $x$ and $y$ is defined by
\be\label{k} k(x,y)=\inf_{\gamma}\{\ell_{k }(\gamma)\},\ee
where the infimum is taken over all rectifiable curves $\gamma$  in $X$ joining $x$ and $y$.

For later use, we need the following elementary estimates in \cite[(2.4)]{BHK}.
\be\label{t-1} k (x,y)\geq \log\Big(1+\frac{d(x,y)}{\min\{d (x), d (y)\}}\Big):=j(x,y)\geq \Big|\log \frac{d (x)}{d (y)}\Big|.\ee

We remark that we do \textbf{not} assume the metric space $X$ to be locally compact. So in general, $X$ equipped with its quasihyperbolic metric is a length space but not geodesic. For this, we need the following notion.

Let $\epsilon\geq 0$.  Then a curve $\alpha$ joining $x$ and $y$ in $X$ is called $\epsilon$-{\it short} if
$$\ell_{k }(\alpha)\leq k (x,y)+\epsilon.$$
It is not hard to see that every subcurve of an $\epsilon$-short curve is also $\epsilon$-short.

\bdefe
Let $A\geq 1$. A rectifiably connected and incomplete metric space $(X,d)$ is called $A$-{\it uniform} provided each pair of points $x$, $y$ in $X$ can
be joined by a rectifiable curve $\gamma$ satisfying
\begin{enumerate}
\item $\ell(\gamma)\leq A\,d(x,y)$, and
\item $\min\{\ell(\gamma[x,z]),\ell(\gamma[z,y])\}\leq A\,d(z)$ for all $z\in \gamma$.
\end{enumerate}
Moreover, $\gamma$ is called an {\it $A$-uniform} curve.
\edefe

 Next, we recall two auxiliary results for our later use.

\begin{lem}\label{Lem-A}$($\cite[Corollary 3.3]{BH}$)$
Let $(X,d)$ be an $A$-uniform metric space. Then for all $x,$ $y$ in $X$, there is a constant $B=B(A)$ such that any $\epsilon$-short curve in $X$ connecting $x$ and $y$ is $B$-uniform provided $\epsilon\leq \epsilon_0$, where $\epsilon_0=\min\{k (x,y),1\}$.
\end{lem}

\begin{lem}\label{Lem-B}
$($\cite[Lemma $2.13$]{BHK}$)$
Let $(X,d)$ be an $A$-uniform space. Then for all $x,$ $y\in X,$ we have
\be\label{z15-1}  k(x,y)\leq 4A^2 \log\Big(1+\frac{d(x,y)}{\min\{d (x), d (y)\}}\Big). \ee
\end{lem}

\subsection{Hyperbolization  metric} In \cite{Z1}, Ibragimov introduced a metric that hyperbolizes $($in the sense of Gromov$)$ all locally compact and incomplete metric spaces, which is roughly quasi-isometric to the quasihyperbolic metric on uniform spaces.

Let $(X,d)$ be a locally complete and incomplete metric space.  Following \cite{Z1, Z2}, we define the {\it hyperbolization  metric} $h$ of $X$ as
\be\label{h} h(x,y)=2\log\frac{d(x,y)+\max\{d(x),d(y)\}}{\sqrt{d(x)d(y)}},\ee
for all $x,y\in X.$

We next introduce the definition of another hyperbolic type metric which is related to the hyperbolization  metric $h$. That is, for all $x,y\in X$, we define
\be\label{j}\widetilde{j}(x,y)=\frac{1}{2}\log \left [ \left (1+\frac{d(x,y)}{d(x)}\right )\left (1+\frac{d(x,y)}{d(y)}\right )   \right ].\ee

Finally, we cite the following auxiliary results.

\begin{thm}\label{zz2}  Let $(X,d)$ be a locally complete and incomplete metric space.
\begin{enumerate}
\item $($\cite[Theorem $2.1$]{Z1}$)$ The space $(X,h)$ is $\delta$-hyperbolic with $\delta=\log 4$.
\item $($\cite[Theorem $3.1$]{Z1}$)$ For all $x,y\in X$, $2\widetilde{j}(x,y)\leq h(x,y)\leq 2\widetilde{j}(x,y)+2\log 2$.
\item $($\cite[Theorem $3.1$]{Z2}$)$ The Gromov boundary $\partial_h X$ of $(X,h)$ is naturally  homeomorphic to its metric boundary $\partial X ~(\partial X \cup \{\infty\}$ if $X$ is unbounded$)$.
\end{enumerate}
\end{thm}

\section{Proofs of Theorem \ref{thm-1} and Theorem \ref{c-1}}\label{sec-3}
\subsection{Proof of Theorem \ref{thm-1}}
Assume that $(X,d)$ is an unbounded locally complete and incomplete metric space, and that $h$ is the hyperbolization metric defined in (\ref{h}). One first observes from Theorem \ref{zz2} that $(X,h)$ is $\delta$-hyperbolic with $\delta=\log4$ and there is a natural bijective map
$$\varphi\colon \partial_h X\to \partial X\cup \{\infty\}$$
with  $\varphi(\xi)=\infty$ for some point $\xi\in \partial_h X$. Let $h_{b,\varepsilon}$ be a Hamenst\"adt metric on the punctured Gromov boundary $\partial_h X\setminus\{\xi\}$ based at the Busemann function $b=b_{w,\xi}\in \mathcal{B}_h(\xi)$, $w\in X$ and parameter $0<\varepsilon\leq \varepsilon_0(\delta)$.

To prove Theorem \ref{thm-1}, it suffices to show that
$$\varphi\colon (\partial_h X\setminus\{\xi\}, h_{b,\varepsilon}) \to (\partial X,d)$$
is $\eta$-quasisymmetric for some function $\eta$ depending only on $\varepsilon$. Since the inverse map of a quasisymmertic map is also quasisymmetric, we only need to check that $\varphi^{-1}$ is $\eta$-quasisymmetric.

To this end, by (\ref{h}), a direct computation guarantees that for all $u,v\in X$,
\beq\label{g-4}\;\;\;\;\;\;\; (u|v)_w &=& \frac{1}{2}\big(h(u,w)+h(v,w)-h(u,v)\big)
\\ \nonumber  &=& \log\frac{[d(w,u)+\max\{d(w),d(u)\}][d(w,v)+\max\{d(w),d(v)\}]}
{d(w)[d(u,v)+\max\{d(u),d(v)\}]}.\eeq

Fix three distinct points $x,y,z\in \partial X$. Take sequences $\{x_n\}, \{y_n\}$, and $ \{z_n\}$ in $X$ converging in the metric $d$ to the points $x,y$ and $z$, respectively. Also, choose a sequence $\{u_n\}$ in $X$ such that $d(u_n,p)\to \infty$ as $n\to\infty$, for all $p\in X$. Thus again by Theorem \ref{zz2}, $\{u_n\}$ is a Gromov sequence in the hyperbolic space $(X,h)$ with $\{u_n\}\in\xi$.  Moreover, we see that $\{x_n\}\in \xi_x=\varphi(x), \{y_n\}\in \xi_y=\varphi(y)$, and $ \{z_n\}\in \xi_z=\varphi(z)$.

Moreover, by (\ref{z0}) and Lemma \ref{z00}, we have for $n$ sufficiently large,

\vspace{8pt}
\noindent
$$(x_n|z_n)_b-(x_n|y_n)_b$$
\begin{align}\label{z-1.1} 
& \doteq_{C(\delta)}    (x_n|z_n)_w-(\xi|z_n)_w-(x_n|y_n)_w+(y_n|\xi)_w   \nonumber\\
&\doteq_{C(\delta)}  (x_n|z_n)_w-(u_n|z_n)_w-(x_n|y_n)_w+(y_n|u_n)_w.
\end{align}

On the other hand, we compute by (\ref{g-4}) that

\vspace{8pt}
\noindent
$$e^{(x_n|z_n)_w-(u_n|z_n)_w-(x_n|y_n)_w+(y_n|u_n)_w}$$
\begin{align}
&=  \frac{d(x_n,y_n)+\max\{d(x_n),d(y_n)\}}{d(x_n,z_n)+\max\{d(x_n),d(z_n)\}}  \cdot\frac{d(u_n,z_n)+\max\{d(u_n),d(z_n)\}}{d(u_n,y_n)+\max\{d(u_n),d(y_n)\}}  \nonumber
\\   &\leq\bigg[\frac{d(x_n,y_n)}{d(x_n,z_n)}+\frac{\max\{d(x_n),d(y_n)\}}{d(x_n,z_n)}\bigg] \cdot\bigg[1+\frac{2d(y_n,z_n)}{d(u_n,y_n)+d(u_n)}\bigg].  \nonumber
\end{align}

This, together with (\ref{z-1.1}) and \cite[Lemma 3.2.4]{BuSc}, shows that there is a constant $C_1=C_1(\delta)$ such that
\be\label{z-1.2} e^{(\xi_x|\xi_z)_b-(\xi_x|\xi_y)_b}\leq C_1 \frac{d(x,y)}{d(x,z)}.\ee
Thus by (\ref{z-0.1}) and (\ref{z-1.2}), we obtain
$$\frac{h_{b,\varepsilon}(\xi_x,\xi_y)}{h_{b,\varepsilon}(\xi_x,\xi_z)}\leq 2e^{\varepsilon[(\xi_x|\xi_z)_b-(\xi_x|\xi_y)_b]}\leq 2C_1^\varepsilon \Big(\frac{d(x,y)}{d(x,z)}\Big)^{\varepsilon}.
$$
Setting $C=2C_1^\varepsilon$ and $\eta(t)=Ct^\varepsilon$. This ensures that $\varphi^{-1}$ is $\eta$-quasisymmetric,  completing the proof.
\qed

\subsection{Tools for the proof of Theorem \ref{c-1}.}
Let $(X,d)$ and $(X',d')$ be metric spaces. Also, we let $\lambda>0 $ and  $\mu\geq 0$. Following \cite{BS}, a map $f\colon X\to X'$ is said to be  a {\it $(\lambda,\mu)$-rough similarity} if every point $x'\in X'$ has distance at most $\mu$ from $f(X)$ and for all $x,y\in X$,
$$\lambda d(x,y)-\mu\leq d'(f(x),f(y))\leq \lambda d(x,y)+\mu.$$

To prove Theorem \ref{c-1}, we need an auxiliary result. Note that the spaces in Lemma~\ref{zz-01} below are not assumed to be geodesic.

\begin{lem}\label{zz-01} Let $\lambda>0$, $\mu,\delta\geq 0$, and let $f\colon X\to X'$ be a $(\lambda,\mu)$-rough similarity from a $\delta$-hyperbolic space $(X,d)$ to a metric space $(X',d')$. Then we have the following:
\begin{enumerate}[(1)]
\item $X'$ is $\delta'$-hyperbolic for some constant $\delta'=\delta'(\lambda,\mu,\delta)\geq 0$.
\item $f$ induces an $\eta$-quasisymmetric bijective map
$$\partial f\colon (\partial_\infty X\setminus\{\xi\},d_{b,\varepsilon})\to (\partial_\infty X'\setminus\{\xi'\},d'_{b',\varepsilon'})$$
with $\eta$ depending only on $\lambda,\mu,\delta,\varepsilon$ and $\varepsilon'$, where $d_{b,\varepsilon}$ and $d'_{b',\varepsilon'}$ are Hamenst\"adt metrics on $\partial_\infty X\setminus\{\xi\}$ and $\partial_\infty X'\setminus\{\xi'\}$, respectively, with $b=b_{\xi,o}$, $b'=b'_{\xi',o'}$, $o\in X$, $o'\in X'$ and $\xi'=\partial f(\xi)$.
\end{enumerate}
\end{lem}

\bpf We first prove (1). For all points $x_1',x_2',x_3',p'\in X'$, since $f$ is $(\lambda,\mu)$-roughly similar, we see that there are points $x_1,x_2,x_3,p\in X$ such that
$$\max_{i=1,2,3} d'(f(x_i),x_i')\leq \mu\;\;\;\;\;\mbox{and}\;\;\;\;\; d'(f(p),p')\leq \mu.$$
Thus a direct computation gives that
\be\label{t-4} (x_i'|x_j')_{p'} \doteq_{C(\mu)} (f(x_i)|f(x_j))_{f(p)} \doteq_{C(\mu)} \lambda (x_i|x_j)_p, \ee
for all $i\neq j\in\{1,2,3\}$. Moreover, since $X$ is $\delta$-hyperbolic, we get by (\ref{t-4}) that
\beq\nonumber (x_1'|x_3')_{p'} &\doteq_{C(\mu)}& \lambda (x_1|x_3)_p
\\ \nonumber&\geq& \lambda \min\{(x_1|x_2)_p,(x_2|x_3)_p \}-\lambda\delta
\\ \nonumber&\doteq_{C(\lambda,\mu)}& \min\{(x_1'|x_2')_{p'},(x_2'|x_3')_{p'} \}-\lambda\delta,
\eeq
 which implies (1).

For (2), we next show that $f$ induces a well-defined bijective map $\partial f\colon \partial_\infty X\to \partial_\infty X'.$
To this end, fix a base point $o\in X$ and for a Gromov sequence $\{x_i\}\subset X$ be given. Since $f$ is $(\lambda,\mu)$-roughly similar,  one observes that
$$(f(x_i)|f(x_j))_{f(o)} \doteq_{C(\mu)} \lambda (x_i|x_j)_o\to \infty\;\;\;\;\mbox{as}\; i,j\to \infty.$$
Thus, we find that the sequence $\{f(x_i)\}$ is also a Gromov sequence in the hyperbolic space $X'$. Moreover, if $\{x_i\}$ and $\{y_i\}$ are equivalent Gromov sequences in $X$, then a similar argument as above shows that $\{f(x_i)\}$ and $\{f(y_i)\}$ are also equivalent.

Therefore, we get a well-defined map, induced by $f$,
$$\partial f\colon \partial_\infty X\to \partial_\infty X'$$
 with $\partial f(\zeta)=\zeta'$ for all $\zeta\in \partial_\infty X$, where $\zeta'\in\partial_\infty X'$ is the equivalence class of $\{f(x_i)\}$ for $\{x_i\}\in \zeta$.

Now we prove that $\partial f$ is injective. If $\partial f(\zeta_1)=\partial f(\zeta_2)$, then for any Gromov sequences $\{x_i\}\in \zeta_1$ and $\{y_i\}\in \zeta_2$,
$$\lambda (x_i|y_i)_o \doteq_{C(\mu)}(f(x_i)|f(y_i))_{f(o)}\to \infty\;\;\;\;\mbox{as}\;i\to \infty,$$
because $f$ is $(\lambda,\mu)$-roughly similar. This implies that $\zeta_1=\zeta_2$, as desired.

We also need to verify that $\partial f$ is surjective. For a given point $\zeta'\in \partial_\infty X'$, take a Gromov sequence $\{x_i'\}\in \zeta'$. Thus for each $i$, there is a point $x_i\in X$ such that $d'(f(x_i),x_i')\leq \mu$. Moreover, since $f$ is $(\lambda,\mu)$-roughly similar, we find that
$$ \lambda (x_i|x_j)_o\doteq_{C(\mu)} (f(x_i)|f(x_j))_{f(o)}\doteq_{C(\mu)}  (x_i'|x_j')_{f(o)}.$$
It follows that both $\{x_i\}$ and $\{f(x_i)\}$ are Gromov sequences with $\{x_i\}\in \zeta$ for some point $\zeta\in\partial_\infty X$. Since
$$(f(x_i)|x_i')_{f(o)}\geq d'(f(o),x_i')-\mu\to \infty \;\;\;\;\mbox{as}\;i\to \infty,
$$
we know that $\{f(x_i)\}\in \zeta'$, which implies that $\partial f(\zeta)=\zeta'$. Hence we  get a well-defined bijective map $\partial f\colon \partial_\infty X\to \partial_\infty X'.$

Finally, it remains to show that the homeomorphism $\partial f\colon (\partial_\infty X\setminus\{\xi\},d_{b,\varepsilon})\to (\partial_\infty X'\setminus\{\xi'\},d'_{b',\varepsilon'})$ is $\eta$-quasisymmetric.

Fix three distinct points $x,y,z\in \partial_\infty X\setminus\{\xi\}$. Choose Gromov sequences $\{x_n\}\in x, \{y_n\}\in y, \{z_n\}\in z$, and $ \{u_n\}\in \xi$. From the definition of $\partial f$, it follows that $\{f(x_n)\}\in \partial f(x)=x', \{f(y_n)\}\in \partial f(y)=y', \{f(z_n)\}\in \partial f(z)=z'$,  and  $ \{f(u_n)\}\in \partial f(\xi)=\xi'$. 

Moreover, by (\ref{z0}) and Lemma \ref{z00}, we have for $n$ sufficiently large,
\begin{align} \label{zz-02}
(x_n|y_n)_b-(x_n|z_n)_b &\doteq_{C(\delta)} (x_n|y_n)_o-(\xi|y_n)_o-(x_n|z_n)_o+(z_n|\xi)_o  \nonumber \\
&\doteq_{C(\delta)} (x_n|y_n)_o-(u_n|y_n)_o-(x_n|z_n)_o+(z_n|u_n)_o.
\end{align}
Since $f$ is $(\lambda,\mu)$-roughly similar and $\{f(u_n)\}\in \partial f(\xi)=\xi'$, again by (\ref{z0}), Lemma~\ref{z00} and (\ref{zz-02}), we have for $n$ sufficiently large,

\vspace{8pt}
\noindent
$$(f(x_n)|f(y_n))_{b'}-(f(x_n)|f(z_n))_{b'}$$
\begin{align}\label{zz-03}
&\doteq_{C(\delta)} \hspace{.2cm} (f(x_n)|f(y_n))_{o'}-(f(y_n)|f(u_n))_{o'}-(f(x_n)|f(z_n))_{o'}+(f(u_n)|f(z_n))_{o'} \nonumber\\
&\;\;\;= \hspace{.5cm} (f(x_n)|f(y_n))_{f(o)}-(f(y_n)|f(u_n))_{f(o)}\nonumber \\
& \hspace{5cm} -(f(x_n)|f(z_n))_{f(o)}+(f(u_n)|f(z_n))_{f(o)} \nonumber\\
&\doteq_{C(\delta,\lambda,\mu)} ~\lambda [(x_n|y_n)_o-(u_n|y_n)_o-(x_n|z_n)_o+(z_n|u_n)_o] \nonumber\\
&\doteq_{C(\delta,\lambda,\mu)} ~\lambda[(x_n|y_n)_b-(x_n|z_n)_b].
\end{align}

Therefore, by (\ref{zz-03}) and \cite[Lemma 3.2.4]{BuSc}, we find that
\be\label{zz-04}(x'|y')_{b'}-(x'|z')_{b'}\doteq_{C(\delta,\lambda,\mu)} \lambda[(x|y)_b-(x|z)_b].\ee

Furthermore, by (\ref{z-0.1}) and (\ref{zz-04}), we obtain
\beq\nonumber \frac{d'_{b',\varepsilon'}(x',y')}{d'_{b',\varepsilon'}(x',z')}&\leq& 2e^{\varepsilon'[(x'|y')_{b'}-(x'|z')_{b'}]}
\\ \nonumber&\leq& C_1 e^{\lambda \varepsilon'[(x|y)_b-(x|z)_b]}
\\ \nonumber&\leq& C_2 \Big[\frac{d_{b,\varepsilon}(x,y)}{d_{b,\varepsilon}(x,z)}\Big]^{\lambda \varepsilon'/\varepsilon}
\eeq
 with constants $C_1$ and $C_2$ depending only on $\lambda,\mu,\delta,\varepsilon$ and $\varepsilon'$.
By setting $\eta(t)=C_2(t)^{\lambda \varepsilon'/\varepsilon}$,  one immediately sees that $\partial f$ is $\eta$-quasisymmetric.  This completes the proof of Lemma \ref{zz-01}.
\epf

 Now we are ready to show Theorem \ref{c-1}.

\subsection{Proof of Theorem \ref{c-1}.}  Let $(X,d)$ be an unbounded locally complete and incomplete metric space, and let $\widetilde{j}$ be defined in \eqref{j}. We first see from Theorem \ref{zz2} that $(X,h)$ is $\delta$-hyperbolic with $\delta=\log 4$ and the identity map $(X,h)\to (X,\widetilde{j})$ is $(2,2\log 2)$-roughly similar. Then it follows from Lemma \ref{zz-01}  that $(X,\widetilde{j})$ is $\delta'$-hyperbolic for some positive number $\delta'$ and the identity map induces an $\eta_1$-quasisymmetric bijective map
$$\psi\colon (\partial_{\widetilde{j}} X\setminus\{\xi\},\widetilde{j}_{b,\varepsilon})\to (\partial_h X\setminus\{\xi'\},h_{b',\varepsilon'})$$
with $\eta_1$ depending only on $\varepsilon$ and $\varepsilon'$, where $\widetilde{j}_{b,\varepsilon}$ and $h_{b',\varepsilon'}$ are Hamenst\"adt metrics on the punctured Gromov boundaries $\partial_{\widetilde{j}} X\setminus\{\xi\}$ and $\partial_h X\setminus\{\xi'\}$ with $b=b_{\xi,o}\in \mathcal{B}_{\widetilde{j}}(\xi)$, $b'=b'_{\xi',o'}\in \mathcal{B}_h(\xi')$, $o\in X$, $o'\in X'$, and $\xi'=\psi(\xi)$, respectively.

On the other hand, one observes from Theorem \ref{thm-1} that there is a natural $\eta_2$-quasisymmetric map
$$\varphi\colon (\partial_h X\setminus\{\xi'\},h_{b',\varepsilon'})\to (\partial X,d)$$
 with $\varphi(\xi')=\infty$ and the function $\eta_2$ depending only on $\varepsilon'$.

Therefore, we obtain a natural $\eta$-quasisymmetric map
$$\widetilde{\varphi}=\varphi\circ \psi\colon (\partial_{\widetilde{j}} X\setminus\{\xi\},\widetilde{j}_{b,\varepsilon})\to (\partial X,d)$$
with  $\widetilde{\varphi}(\xi)=\infty$ and $\eta=\eta_2\circ \eta_1$, because the composition of quasisymmetric maps are also quasisymmetric.
\qed

\section{Proofs of Theorem \ref{thm-2} and Theorem \ref{c-2}}\label{sec-4}
\subsection{Hyperbolizing unbounded uniform spaces}
In this subsection, we assume that $(X,d)$ is an unbounded $A$-uniform metric space for some constant $A\geq1$ and $k$ is the quasihyperbolic metric of $X$. We begin with a non-locally compact version of \cite[Theorem 3.6]{BHK} and a metric version of \cite[Theorem 2.12]{Vai}, i.e., the Gromov hyperbolicity of a uniform space with respect to its quasihyperbolic metric.

By Lemma \ref{Lem-A}, we see that for every pair of points in a uniform space, there is a $1$-short curve $\gamma$ joining them such that $\gamma$ is a uniform curve, which, in what follows, is briefly called a {\it $1$-short uniform} curve. Thus, by replacing ``$2$-neargeodesic" by ``$1$-short uniform curve", the similar reasoning as in the proof of \cite[Theorem 2.12]{Vai} shows that the following result holds true.

\begin{lem}\label{thm-3-Lem} Let $(X,d)$ be an $A$-uniform metric space and $k$ its quasihyperbolic metric. Then $(X,k)$ is $\delta$-hyperbolic with $\delta=\delta(A)$.
\end{lem}

As the space $(X,d)$ is not assumed to be locally compact, $(X,k)$ is a length metric space but not geodesic. Thus,  we need the following {\it standard estimate} established by V\"ais\"al\"a for length Gromov hyperbolic spaces.

\begin{lem}$($\cite[2.33]{Vai-0}$)$
Suppose that $(X,k )$ is $\delta$-hyperbolic and that $\alpha$ is a quasihyperbolic $\epsilon$-short curve connecting $x$ to $y$.
Then for all $p\in X$,
\be\label{zz-05} k (p,\alpha)-2\delta-\epsilon\leq (x|y)_p\leq k (p,\alpha)+\frac{\epsilon}{2}. \ee
Note that here,
$$(x|y)_p=\frac{1}{2}[k(x,p)+k(y,p)-k(x,y)].$$
\end{lem}

Moreover, the metric version of \cite[Lemma 2.22]{Vai} is needed.

\begin{lem}\label{zz-06} The natural map $\varphi\colon  \partial_k X \rightarrow \partial{X}  \cup\{\infty\}$ exists if and only if every Gromov sequence in $X$ has a limit in metric topology, where $\partial_k X$ is the Gromov boundary of hyperbolic space $(X,k)$.
\end{lem}

The following result is also a metric version of \cite[Lemma $2.25$]{Vai}. The proof follows  on the same line  but for the sake of completeness we include the details.

\begin{lem}\label{zz-07} Let $(X,d)$ be an unbounded $A$-uniform metric space and $k$ its quasihyperbolic metric. For all $p\in X$ and $u\in \partial X\cup\{\infty\}$, we have $(x|y)_p\rightarrow \infty$ as $x,y\rightarrow u$ in the metric topology.
\end{lem}

\bpf We first assume that $u\in \partial X$. Let $r>0$ and $x,$ $y\in B(u,r)$. By Lemma \ref{Lem-A}, there is a $1$-short $B$-uniform curve $\alpha$ joining $x$ and $y$ with $B=B(A)$. Fix $z\in \alpha$.  To show that $(x|y)_p\to \infty$ as $x, y\to u$, we know from the standard estimate \eqref{zz-05} that it suffices to find an estimate $k (p,z)\geq M(r)$ such that
 $$M(r)\rightarrow \infty\;\;{\rm as}\;\;r\rightarrow 0.$$

To this end, since $\alpha$ is $B$-uniform and $u\in\partial X$, we compute
$$d(z)\leq d(z,u)\leq d(z,x)+d(x,u)\leq Bd(x,y)+r\leq (2B+1)r.$$
This implies that
$$k (p,z)\geq \log\frac{d (p)}{d (z)}\geq \log\frac{d (p)}{(2B+1)r}=M(r),$$
as required.

We are thus left to consider the case that $u=\infty$. Choose a point $v\in \partial X$ and let $R>B d(v,p)$. For $x,y\in X\setminus B(p,R)$,  Lemma \ref{Lem-A}  ensures that there is a $1$-short $B$-uniform curve $\beta$ connecting $x$ to $y$. Fix $z\in \beta$.  To prove that $(x|y)_p\to \infty$ as $x, y\to \infty$, we only need to show that
$$k (z,p)\geq M(R)\rightarrow \infty\;\;{\rm  as}\;\; R\rightarrow \infty.$$
Since $\beta$ is $B$-uniform, we get
$$R-d(z,p)\leq \min\{d(x,z), d(z,y)\} \leq Bd (z)\leq Bd(z,v)\leq Bd(v,p)+Bd(p,z),$$
and therefore, $$d(z,p)\geq \frac{R-Bd(v,p)}{1+B}.$$
This guarantees that
$$k (z,p)\geq \log\Big(1+\frac{d(p,z)}{d (p)}\Big)\geq \log\Big(1+\frac{R-Bd(v,p)}{(1+B)d (p)}\Big)=M(R),$$
as desired. Hence the lemma follows.
\epf

To prove Theorem \ref{thm-2}, we also need an auxiliary lemma. The proof is the same as \cite[Lemma 2.30]{Vai}.

\begin{lem}\label{zz-09} Let $(X,d)$ be an unbounded $A$-uniform metric space and $x, y, u\in X$ with $d(x,y)\leq \frac{1}{2}d(x,u)$. Let $\alpha$ and $\beta$ be $1$-short $B$-uniform curves connecting $x$ to $u$ and $y$, respectively. Define a point $v$ in $\alpha$ by
$$\ell(\alpha[x,v])=d(x,y).$$
Then we have
$$ k(v,u) \doteq_{C(A)} k (u,\beta).$$
\end{lem}

 Now, we are ready to show Theorem \ref{thm-2}.

\subsection{Proof of Theorem \ref{thm-2}}  Let $A\geq 1$ and let $(X,d)$ be an unbounded $A$-uniform metric space and $k$ its quasihyperbolic metric. By Lemma \ref{thm-3-Lem}, we know that $(X,k)$ is $\delta$-hyperbolic for some constant $\delta=\delta(A)\geq 0$. We first prove the existence of the natural map
$$ \phi\colon \partial_k X\to \partial X\cup \{\infty\}.$$

To this end, let $p\in X$ and let $\{x_i\}$ be a Gromov sequence in the $\delta$-hyperbolic space $(X,k)$. By Lemma \ref{zz-06}, we only need to show that $\{x_i\}$ has a limit in the $d$-metric topology.  We consider two cases based on the boundedness of this sequence.

\vskip 2mm
\noindent
{\bf Case A:} Suppose that $\{x_i\}$ is bounded.
\vskip 2mm

In this case, it suffices to verify that $\{x_i\}$ is a Cauchy sequence in $(X,d)$. For any $i\not=j$, set $d(x_i,x_j):=t_{ij}$. By Lemma \ref{Lem-A}, we know that there is a $1$-short $B$-uniform curve $\gamma_{i,j}$ joining $x_i$ and $ x_j$. Now take some point $z$ in $\gamma_{i,j}$ satisfying
\be\label{t-2} \ell(\gamma[x_i,z])=\frac{1}{2}t_{ij}.\ee
Since $\gamma_{i,j}$ is $B$-uniform, we have
$$d (z)\geq \frac{1}{B}\ell(\gamma[x_i,z])= \frac{1}{2B}t_{ij}.$$

Moreover, since $\{x_i\}$ is bounded, there is a positive number $M>0$ such that
$d( x_i,p)\leq M$ for all positive integer $i$.
Thus by (\ref{t-2}), we compute
$$d(z,p)\leq d(z,x_i)+d(x_i,p)\leq \frac{1}{2}(d(x_i,p)+d(x_j,p))+d(x_i,p)\leq 2M.$$
 Combining with \eqref{z15-1} and \eqref{zz-05}, we get
\begin{eqnarray*}(x_i|x_j)_p&\leq& k (p,z)+\frac{1}{2}
\\ &\leq& 4A^2\log\Big(1+\frac{d(p,z)}{\min\{d (p),d (z)\}}\Big)+\frac{1}{2}
\\ &\leq& 4A^2\log\Big(1+\frac{2M}{\min\{d (p),t_{ij}/2B\}}\Big)+\frac{1}{2}.
\end{eqnarray*}
Since $(x_i|x_j)_p\rightarrow \infty$ as $i,j\rightarrow \infty$,  it follows that $t_{ij}\rightarrow 0$, which shows that $\{x_i\}$ is a $d$-Cauchy sequence.

\vskip 2mm
\noindent
{\bf Case B:} Suppose that $\{x_i\}$ is unbounded.
\vskip 2mm

In this case, we show that $d(x_i,p)\rightarrow \infty$ as $i\to \infty$. Suppose on the contrary that there is some $R>0$ such that $d(x_i,p)\leq R$ for infinitely many $i$. Let $x_s,x_t\in\{x_i\}$ be such that $d(x_s,p)\leq R$ and $d(x_t,p)\geq 3R$.  Again by Lemma \ref{Lem-A}, there is a $1$-short $B$-uniform curve $\beta$ joining $x_s$ and $x_t$. Choose a point $z\in \beta$ with $d(z,p)=2R$. Since $\beta$ is $B$-uniform, we have
$$d(z)\geq \frac{1}{B}\min\{\ell(\beta[x_s,z]), \ell(\beta[x_t,z])\}\geq \frac{R}{B}.$$
This, together with \eqref{z15-1}, gives that
$$k (p,\beta)\leq k (p,z)\leq 4A^2\log\Big(1+\frac{2R}{\min\{d (p),R/B\}}\Big):=R_0.$$
Again by \eqref{zz-05}, we get
$$(x_s|x_t)_p\leq R_0+\frac{1}{2},$$
which is abused because $\{x_i\}$ is a Gromov sequence. Hence, the existence of the natural map $\phi$ follows.

Next, we check the injectivity of  $\phi$. Towards this end, we let $ \xi, \zeta\in\partial_k X$ with $\phi(\xi)=\phi(\zeta)=u$ and let $\{x_i\}\in \xi$ and  $\{y_j\}\in \zeta$ be Gromov sequences. By Lemma~\ref{zz-07}, we know that $\{x_i\}$ is equivalent to $\{y_j\}$. Hence $\xi=\zeta$.

We also need to verify that  $\phi$ is surjective. Let $u\in\partial X$ and choose a sequence  $\{x_i\}$  in $X$ converging to $u$ in metric topology. Again by Lemma \ref{zz-07}, we have $(x_i|x_j)_p\rightarrow \infty$ as $i,j\rightarrow \infty$. Thus   $\{x_i\}$  is a Gromov sequence  with $\{x_i\}\in \xi\in \partial_{\infty}X$ and $\phi(\xi)=u$. It follows that  $\phi$ is surjective.

Because the inverse map of a quasisymmetric mapping is also quasisymmetric, we only need to show that  $\phi^{-1}$ is $\eta$-quasisymmetric. For any distinct points $x,y,z\in \partial X$ with $d(x,y)=t d(x,z)$, it remains to prove that
$$k_{b,\varepsilon}(x,y)\leq \eta(t)k_{b,\varepsilon}(x,z)$$
for some homeomorphism $\eta\colon [0,\infty)\to [0,\infty)$, where $k_{b,\varepsilon}$ is a Hamenst\"adt metric on the punctured Gromov boundary $\partial_k X\setminus\{\xi\}$ based at a Busemann function $b=b_{\xi,o}\in \mathcal{B}_k(\xi)$ with $\phi(\xi)=\infty$ and parameter $0<\varepsilon\leq \varepsilon_0(\delta)$.

Fix three distinct points $x,y,z\in \partial X$. Take sequences $\{x_n\}, \{y_n\}$, and $ \{z_n\}$ in $X$ converging in the metric $d$ to the points $x, y$, and $z$, respectively. Also, choose a sequence $\{u_n\}$ in $X$ such that $d(u_n,p)\to \infty$ as $n\to\infty$, for all $p\in X$. By the argument in \textbf{Case B}, we know that $\{u_n\}$ is a Gromov sequence with $\{u_n\}\in\xi$. Similarly, we see that $\{x_n\}\in \xi_x=\phi^{-1}(x), \{y_n\}\in \xi_y=\phi^{-1}(y)$, and $\{z_n\}\in \xi_z=\phi^{-1}(z)$.
Using (\ref{z0}), we have
\be\label{z-1} (x_n|z_n)_b-(x_n|y_n)_b\doteq_{C(\delta)} (x_n|z_n)_o-(\xi|z_n)_o-(x_n|y_n)_o+(y_n|\xi)_o. \ee
Moreover, by Lemma \ref{z00}, we have for $n$ sufficiently large,

\vspace{8pt}
\noindent
$(x_n|z_n)_o-(\xi|z_n)_o-(x_n|y_n)_o+(y_n|\xi)_o
$
\begin{align} \label{z-2}
&\doteq_{C(\delta)} (x_n|z_n)_o-(u_n|z_n)_o-(x_n|y_n)_o+(y_n|u_n)_o \nonumber\\
&= (x_n|z_n)_{u_n}-(x_n|y_n)_{u_n}.
\end{align}
Now by (\ref{z-1}) and (\ref{z-2}), for $n$ sufficiently large, we see that there is a constant $C_1=C_1(\delta)$ such that
\be\label{z-3} T_n:=\frac{e^{-\varepsilon(x_n|y_n)_b}}{e^{-\varepsilon(x_n|z_n)_b}}
=e^{\varepsilon[(x_n|z_n)_b-(x_n|y_n)_b]}\leq C_1 e^{\varepsilon[(x_n|z_n)_{u_n}-(x_n|y_n)_{u_n}]}.
\ee

Fix a sufficiently large natural number $n$, we let
$$s_n=(x_n|z_n)_{u_n}-(x_n|y_n)_{u_n}\;\;\;\;\mbox{and} \;\;\;\;t_n= \frac{d(x_n,y_n)}{d(x_n,z_n)}.$$
Since the sequence $\{u_n\}$ $d$-converges to $\infty$, we may assume without loss of generality that
\be\label{z-4} \max\{d(x_n,y_n), d(x_n,z_n)\}\leq \frac{1}{2} d(x_n,u_n).\ee
Furthermore, by Lemma \ref{Lem-A}, we find that there are $1$-short $B$-uniform curves $\alpha_n, \beta_n$,  and $\gamma_n$ connecting $x_n$ to the points $u_n, y_n$, and $z_n$, respectively. By (\ref{z-4}), we may pick two points  $p_n$ and $q_n$ in $\alpha_n$ such that
\be\label{z-5} \ell(\alpha[x_n,p_n])=d(x_n,y_n)\;\;\;\;\;\mbox{and} \;\;\;\;\;\ell(\alpha[x_n,q_n])=d(x_n,z_n).\ee
As $\alpha_n$ is $B$-uniform, we get
\be\label{z-6} d(p_n)\geq \frac{1}{B} d(x_n,y_n) \;\;\;\;\;\;\mbox{and} \;\;\;\;d(q_n)\geq \frac{1}{B} d(x_n,z_n).\ee

In the following, we consider two cases.

\vskip 2mm
\noindent
{\bf Case I:} Suppose that $t\geq 1$.
\vskip 2mm

In this case, we may assume $t_n\geq 1/2$ for $n$ sufficiently large, because $t_n\to t$ as $n\to \infty$. Thus we have
$$d(x_n,y_n)\geq \frac{1}{2} d(x_n,z_n),$$
from which and (\ref{z-6}) it follows that
\be\label{z-7} \min\{d(p_n),d(q_n)\}\geq \frac{1}{2B} d(x_n,z_n).\ee
Moreover, by (\ref{z-5}) we get
$$d(p_n,q_n)\leq d(p_n,x_n)+d(x_n,q_n)\leq d(x_n,y_n)+d(z_n,x_n)\leq 3d(x_n,y_n).$$
Thus we deduce from (\ref{z-7}) and Lemma \ref{Lem-B} that
$$k(p_n,q_n)\leq 4A^2\log\Big(1+\frac{d(p_n,q_n)}{\min\{d(p_n),d(q_n)\}}\Big)\leq 4A^2\log(1+6B t_n).$$

Therefore, by (\ref{zz-05}) and Lemma \ref{zz-09}, we obtain
\beq\nonumber s_n &=& (x_n|z_n)_{u_n}-(x_n|y_n)_{u_n}
\\ \nonumber&\doteq_{C(\delta)}& k(u_n,\gamma_n)-k(u_n,\beta_n)
\\ \nonumber&\doteq_{C(A)}& k(u_n,q_n)-k(u_n,p_n)
\\ \nonumber&\leq& k(p_n,q_n)\leq 4A^2\log(1+6B t_n).
\eeq
This, together with (\ref{z-3}), shows that there is a constant $C_2=C_2(A,\varepsilon)$ such that
$$T_n\leq C_2 e^{4\varepsilon A^2\log(1+6B t_n)}.$$
By taking $n\to \infty$, we see from (\ref{z-0.1}) and \cite[Lemma 3.2.4]{BuSc}  that there is a constant $C=C(A,\varepsilon)$ such that
$$\frac{k_{b,\varepsilon}(\xi_x,\xi_y)}{k_{b,\varepsilon}(\xi_x,\xi_z)}\leq Ct^{4A^2\varepsilon},$$
as desired.

\vskip 2mm
\noindent
{\bf Case II:} Suppose that $0<t< 1$.
\vskip 2mm

Since $0<t< 1$ and $t_n\to t$ as $t\ra \infty$, we may assume  without loss of generality that $t_n\leq 1$ for $n$ sufficiently large. Thus by our choices of $p_n$ and $q_n$ in the curve $\alpha_n$, we find that $$q_n\in\alpha_n[p_n,u_n].$$
Moreover, since $\alpha_n$ is $1$-short, again by (\ref{zz-05}) and Lemma \ref{zz-09}, we have
\be\label{z-8}  s_n \doteq_{C(A)} k(u_n,q_n)-k(u_n,p_n)\doteq -k(p_n,q_n).\ee

On the other hand, by (\ref{z-6}) and (\ref{z-5}), we compute
$$d(p_n)\leq d(p_n,x)\leq d(p_n,x_n)+d(x_n,x)\leq d(x_n,y_n)+d(x_n,x).$$
Since $d(x_n,x)\to 0$ as $n\to \infty$, for $n$ sufficiently large we may assume $d(x_n,x)\leq d(x_n,y_n)$. Therefore, we get
$$-k(p_n,q_n)\leq \log \frac{d(p_n)}{d(q_n)} \leq \log\frac{2Bd(x_n,y_n)}{d(x_n,z_n)}=\log(2Bt_n).$$
It follows from  (\ref{z-3}) and (\ref{z-8}) that there is a constant $C_3=C_3(A,\varepsilon)$ such that
$$T_n\leq C_3 e^{\varepsilon \log (2Bt_n)}.$$
Letting $n\to \infty$, again by (\ref{z-0.1}) and \cite[Lemma 3.2.4]{BuSc}, we find that there is a constant $C=C(A,\varepsilon)$ such that
$$\frac{k_{b,\varepsilon}(\xi_x,\xi_y)}{k_{b,\varepsilon}(\xi_x,\xi_z)}\leq Ct^{\varepsilon},$$
as required.
\qed

\subsection{Proof of Theorem \ref{c-2}} Because the necessity follows from Theorem \ref{thm-2}, we only need to show the sufficiency. Assume that $(D,k)$ is $\delta$-hyperbolic and there is a natural $\eta$-quasisymmetric map
$$\phi\colon (\partial_k D\setminus\{\xi\},k_{b,\varepsilon})\to (\partial D,|\cdot|),$$
where $\partial_k D\setminus\{\xi\}$ is the punctured Gromov boundary  of the hyperbolic space $(D,k)$ equipped with a Hamenst\"adt metric $k_{b,\varepsilon}$ and $b=b_{\xi,o}\in \mathcal{B}_k(\xi)$ is a Busemann function with  $\phi(\xi)=\infty$. We check the uniformity condition of $D$.

Firstly, by \cite[Theorem 3.27]{Vai}, it suffices to prove that there is a natural $\theta$-quasim\"obius bijection
$$\widetilde{\phi}\colon (\partial_k D\setminus \{\xi\},k_{p,\epsilon})\to (\partial D, |\cdot|),$$
where $k_{p,\epsilon}$ is a visual metric on $\partial_k D$ with base point $p\in D$ and parameter $\epsilon\leq \min\{1,\frac{1}{5\delta}\}$.

The existence of such a natural map $\phi\colon \partial_k D\to \partial D\cup \{\infty\}$ follows from the assumption and what remains to show is that the identity map
$$\id\colon (\partial_k D\setminus \{\xi\},k_{p,\epsilon})\to (\partial_k D\setminus\{\xi\},k_{b,\varepsilon})$$
is $\theta_0$-quasim\"obius, because we have assumed that $\phi\colon (\partial_k D\setminus\{\xi\},k_{b,\varepsilon})\to (\partial D,|\cdot|)$ is $\eta$-quasisymmetric and the composition of quasisymmetric and quasim\"obius maps are also quasim\"obius.

For any distinct points $x, y, z, w\in \partial_k D\setminus \{\xi\}$, we denote the cross ratio of these four points associated with the metric $k_{p,\epsilon}$ by
$$t=|x,y,z,w|_{k_{p,\epsilon}}.$$
Choose Gromov sequences $\{x_n\}\in x, \{y_n\}\in y, \{z_n\}\in z,$ and $\{w_n\}\in w$. By (\ref{z0}), we compute

\vspace{8pt}
\noindent
$$ (x_n|y_n)_b+(z_n|w_n)_b-(x_n|z_n)_b-(w_n|y_n)_b$$
\begin{align}
&\doteq_{C(\delta)} (x_n|y_n)_o+(z_n|w_n)_o-(x_n|z_n)_o-(w_n|y_n)_o \nonumber\\
&= (x_n|y_n)_p+(z_n|w_n)_p-(x_n|z_n)_p-(w_n|y_n)_p.\nonumber
\end{align}
Thus by \cite[Lemma 3.2.4]{BuSc} and Lemma \ref{z00}, it follows that
$$(x|y)_b+(z|w)_b-(x|z)_b-(w|y)_b \doteq_{C(\delta)}(x|y)_p+(z|w)_p-(x|z)_p-(w|y)_p.$$
Therefore, by (\ref{t-3}) and (\ref{z-0.1}), we find that there is a constant $C=C(\delta)$ such that
\beq\nonumber |x,y,z,w|_{k_{b,\varepsilon}}&=& \frac{k_{b,\varepsilon}(x,z)k_{b,\varepsilon}(y,w)}{k_{b,\varepsilon}(x,y)k_{b,\varepsilon}(z,w)}
\\ \nonumber&\leq& 4e^{\varepsilon[(x|y)_b+(z|w)_b-(x|z)_b-(w|y)_b]}
\\ \nonumber&\leq& Ce^{\varepsilon[(x|y)_p+(z|w)_p-(x|z)_p-(w|y)_p]}
\\ \nonumber&\leq& 4C\Big[\frac{k_{p,\epsilon}(x,z)k_{p,\epsilon}(y,w)}{k_{p,\epsilon}(x,y)k_{p,\epsilon}(z,w)}\Big]^{\varepsilon/\epsilon}
=4Ct^{\varepsilon/\epsilon}.
\eeq
Put $\theta_0(t)=4Ct^{\varepsilon/\epsilon}$. The proof is complete.
\qed


\end{document}